\documentclass[a4 paper]{article}
\usepackage{fancyhdr}
\usepackage{latexsym}
\usepackage{mathrsfs}
\usepackage{cancel}
\usepackage{color}
\usepackage{multirow}
\usepackage{eso-pic}
\usepackage{dsfont}
\usepackage{amssymb}
\usepackage{graphicx}
\usepackage{amsmath}
\usepackage{setspace}
\usepackage{geometry}
\usepackage{amsthm}
\usepackage{hyperref}

\usepackage{amssymb,amsfonts}
\usepackage[all,arc]{xy}
\usepackage{enumerate}
\usepackage{mathrsfs}

\makeatletter
\def\cleardoublepage{\clearpage\if@twoside \ifodd\c@page\else
   \hbox{}\thispagestyle{empty}\newpage\addtocounter{page}{-1}
   \if@twocolumn\hbox{}\newpage\fi\fi\fi}
\makeatother

\usepackage{graphicx}
\DeclareSymbolFont{lettersA}{U}{txmia}{m}{it}
\DeclareMathSymbol{\pi}{\mathord}{lettersA}{25}
\DeclareMathSymbol{\muup}{\mathord}{lettersA}{22}

\begin{document}
\title{On Behaviors of Maximal Dimension}
\author{Wenjun Niu\footnote{This material is based on work at 2017 Cornell University Math REU Program. The author want to express his acknowledgement to Professor R. Keith Dennis and graduate mentor Ravi Fernando for their guidance and encouragement.}\\ Fudan University \\ Shanghai, China, 200433}
\maketitle

\begin{abstract}
In this paper, we investigate behaviors of Maximal Dimension, a group invariant involving certain configuration of maximal subgroups, which we denote by $\mathrm{MaxDim}$. We prove that in some special cases, $\mathrm{MaxDim}(G\times H)=\mathrm{MaxDim}(G)+\mathrm{MaxDim}(H)$. We also prove a conjecture stated by Ellie Thieu which shows that groups with $m<\mathrm{MaxDim}$ can be constructed from groups with $m<i$.
\end{abstract}

\section{Introduction and background}\label{1}

For a group $G$, a subset $s\subseteq G$ of elements in $G$ is a generating set if $\langle g, g\in s\rangle=G$. $s$ is said to be irredundant if $\langle g,g\in s\setminus\{h\}\rangle\ne \langle g,g\in s\rangle$ for each $h\in s$. Thus for a finite group $G$, any generating set $s$ contains an irredundant generating set, simply by removing the redundant elements one at a time from $s$. Note that this is generally not true for infinite groups, as has been mentioned in \cite{fernando2015inequality} by R. Fernando. In what follows, we consider exclusively only finite groups, which will not be mentioned again in the statements. With these in mind, we introduce two different notions of ``dimension" for a group, whose associations has been studied for long. Denote by $m(G)$ the maximal size of an irredundant generating set of $G$, and $i(G)$ the maximal size of an irredundant set of $G$. By definition one clearly has $m(G)\leq i(G)$. It is not hard to see that these functions are not always equal, for example, if $G=\mathbb{Z}_p\wr \mathbb{Z}_p$, then $m(G)=2$ but $i(G)=p$. 

There are several reasons these functions are called ``dimension". First, one may observe that these functions generalizes dimension of a vector space, on which these conceptions coincide. Second, these functions clearly measure how ``large" a group may be. Finally, thanks to the work of several mathematicians, it has been shown that these functions  behave nicely on groups, e.g., $m(G\times H)=m(G)+m(H)$ and $i(G\times H)=i(G)+i(H)$\footnote{The conclusion for $i$ is more straightforward from Whiston's Lemma. For $m$, see \cite{lucchini2013largest} or \cite{genseqfingrp}.}. 

There is another counterpart of dimension for finite groups, which comes somehow naturally from generating sequences, and will be our main gradient. In order to define this function, let us introduce the following definition:

\newtheorem{DefofGen}{Definition}[section]

\begin{DefofGen}\label{DefGen}
Let $G$ be a group, and $\{H_i\mid 1\leq i \leq n\}$ a set of subgroups of $G$. $\{H_i\mid 1\leq i \leq n\}$ is said to be in \textbf{general position} if for any $1\leq j\leq n$, $\bigcap\limits_{i\ne j} H_i\supsetneq \bigcap\limits_i H_i$.
\end{DefofGen}

In other words, a set of subgroups of $G$ is in general position if the intersection of subgroups from any proper subset strictly contains the intersection of the subgroups from the whole set. Now we define:

\newtheorem{DefofMD}[DefofGen]{Definition}

\begin{DefofMD}\label{DefMD}
Let $G$ be a group, the \textbf{maximal dimension} of $G$, denoted by $\mathrm{MaxDim}(G)$, is defined by:
\begin{equation}
\mathrm{MaxDim}(G)=\max_\mathcal{S} |\mathcal{S}|
\end{equation}
where $\mathcal{S}$ ranges over collections of maximal subgroups of $G$ that are in general position. 
\end{DefofMD}

One may recognize the counterpart of this definition in linear algebra. This function is related to generation of finite groups by the following proposition:

\newtheorem{RelofmM}[DefofGen]{Proposition}

\begin{RelofmM}\label{RelofmM}
Let $G$ be a group, then $m(G)\leq \mathrm{MaxDim}(G)\leq i(G)$.
\end{RelofmM}

\begin{proof}
First, we prove that $m(G)\leq \mathrm{MaxDim}(G)$. Let $s=\{g_1,\ldots, g_n\}$ be an irredundant generating set of $G$ with $n=m(G)$. For any $1\leq i\leq n$, choose $M_i$ maximal such that $s\setminus \{g_i\}\subseteq M_i$, which is possible because $s\setminus \{g_i\}$ generates a proper subgroup. One observe that $g_i\notin M_i$ otherwise $M_i$ would contain a generating set. We claim that $\{M_i,1\leq i\leq n\}$ are in general position. For any $j$, we indeed have $\cap_{i\ne j}M_i\supsetneq \cap_{i}M_i$ because the first subgroup contains $g_j$ while the second doesn't. Hence the claim holds and $m(G)=n\leq \mathrm{MaxDim}(G)$.

Next, we prove that $\mathrm{MaxDim}(G)\leq i(G)$. Let $\mathcal{S}=\{M_1,\ldots, M_k\}$ be maximal subgroups of $G$ with $k=\mathrm{MaxDim}(G)$. By definition, for any $i$, $\cap_{j\ne i}M_j\setminus M_i$ is not empty, hence one can choose $g_i$ in this set. We claim that $s=\{g_1,\ldots, g_m\}$ is irredundant. Indeed, the subgroup generated by $s\setminus \{g_i\}$ is contained in $M_i$ which does not contain $g_i$, hence itself cannot contain $g_i$, so the claim is true. Thus $\mathrm{MaxDim}(G)=k\leq i(G)$.
\end{proof}

\newtheorem{Rem1}[DefofGen]{Remark}

\begin{Rem1}
We say $s$ in the above proof certifies $\mathcal{S}$.
\end{Rem1}

In light of Proposition \ref{RelofmM}, R. Keith Dennis suggested a way of seeking generating sets by looking for maximal subgroups in general position attaining maximal length, and then use the above machinery. But unfortunately, the behavior of this function has not been studied very profoundly and adapting this machinery can be risky. R. Fernando(\cite{fernando2015inequality}), E. Detomi and A. Lucchini(\cite{detomi2016maximal}) exhibit groups with $m\ll \mathrm{MaxDim}$, and proved that $m=\mathrm{MaxDim}$ for a large family of groups. Computation by R. Keith Dennis reveals that this function may have nice behavior on product of groups. We will present a partial result in this direction. But before going to next section, let us introduce the following definition:

\newtheorem{Pht}[DefofGen]{Definition}

\begin{Pht}
Let $G$ be a group. The Frattini Subgroup of $G$, denoted by $\Phi(G)$, is defined as the intersection of all maximal subgroups of $G$.
\end{Pht} 

The importance of $\Phi(G)$ for us is that it is negligible in computing $m$ and $\mathrm{MaxDim}$. We say $G$ is \textbf{Frattini Free} if $\Phi(G)$ is trivial. It is not restrictive to always assume that $G$ is Frattini Free when computing these functions because $G/\Phi(G)$ is always Frattini Free. We will do so when necessary. 

\section{On additivity of $\mathrm{MaxDim}$}\label{2}

In this section, we will prove that $\mathrm{MaxDim}(H\times K)=\mathrm{MaxDim}(H)+\mathrm{MaxDim}(K)$ for appropriate $H$ and $K$, and will present the obstacle to generality. First, we will introduce the following:

\newtheorem{LemofGou}{Lemma}[section]

\begin{LemofGou}[Goursat]\label{LemofGou}
Let $H$ and $K$ be groups. Then maximal subgroups of $H\times K$ come from one of the following two types:

(1). $M\times K$ for $M$ maximal in $H$ or $H\times M'$ for $M'$ maximal in $K$.

(2). There exists a tuple $(S,N,N',\alpha)$, where $S$ is a simple group, $N$ and $N'$ are normal subgroups of $H$ and $K$ respectively, and $\alpha\!:\! H/N\to K/N'\cong S$ is an isomorphism, such that $\Delta\subseteq H\times K$ is defined by $\Delta=\{(h,k)\mid \alpha(\bar{h})=\bar{k}\}$.

\end{LemofGou} 

\newtheorem{Rem2}[LemofGou]{Remark}

\begin{Rem2}
Subgroups in (1) will be called standard subgroups, and subgroups in (2) will be called pullback subgroups\footnote{This terminology comes from the fact that such groups are pullbacks in the category of groups.}. 
\end{Rem2}

Lemma \ref{LemofGou} partially explains the reason why $\mathrm{MaxDim}$ is difficult to understand: one has to deal with pullback subgroups, which involves common simple quotients of $H$ and $K$. One easy observation will be that if pullback subgroups do not exist, in other words, if $H$ and $K$ are coprime, then $\mathrm{MaxDim}(H\times K)=\mathrm{MaxDim}(H)+\mathrm{MaxDim}(K)$. This in fact can be generalized slightly to the following:

\newtheorem{THM1}[LemofGou]{Theorem}
\begin{THM1}\label{THM1}
Let $H$ and $K$ be groups admitting no common nonabelian simple quotient. The following identity holds: 
\begin{center}
{ $\displaystyle \mathrm{MaxDim}(H\times K)=\mathrm{MaxDim}(H)+\mathrm{MaxDim}(K)$
}
\end{center}
\end{THM1}
\noindent with which a corollary:

\newtheorem{CorSol}[LemofGou]{Corollary}
\begin{CorSol}\label{CorSol}
Let $H$ and $K$ be groups such that one of them is solvable, then $\mathrm{MaxDim}(H\times K)=\mathrm{MaxDim}(H)+\mathrm{MaxDim}(K)$.

\end{CorSol}

Thus with the help of Theorem \ref{THM1}, the additivity of $\mathrm{MaxDim}$ on product of finite solvable groups has been settled. We proceed to the proof of Theorem \ref{THM1}, beginning with a lemma:

\newtheorem{CruLem}[LemofGou]{Lemma}
\begin{CruLem}\label{CruLem}
Let $H$ be a group and $N$ be a normal subgroup of $H$. Let $M_1,\ldots,M_k$ be maximal subgroups of $H$ in general position, then it is possible to rearrange them, so that $M_1\cap N,\ldots,M_l\cap N$ are in general position and for any $j>l$, $M_j\cap N$ contains $\cap_{1\leq i\leq l}M_i\cap N$. Moreover, if we let $\pi:H\to H/N$ be the projection and $R=\cap_{1\leq i\leq l}M_i$, then $\pi(R\cap M_j)$ for $j>l$ are in general position in $H/N$.
\end{CruLem}

\begin{proof}
Without loss of generality we assume that $M_i\cap N$ are in general position for $1\leq i\leq l$ and $M_j\cap N$ contains the intersection of them for any $j>l$. So it remains to show that $\pi(R\cap M_j)$ for $j>l$ are in general position. For this purpose, one need only show that for any $l<j\leq k$, $\pi(R\cap M_j)$ does not contain $$\bigcap\limits_{\substack{l<t\leq k\\ t\ne j}}\pi(R\cap M_t)$$ Without loss of generality assume $j=k$. Since clearly $$\bigcap\limits_{l<t\leq k-1}\pi(R\cap M_t)\supseteq \pi(\bigcap\limits_{l<t\leq k-1}M_t\cap R )$$ we need only show that $\pi(R\cap M_k)$ does not contain the latter group, which is equivalent to that $$(R\cap M_k)N\nsupseteq \bigcap\limits_{l<t\leq k-1} M_t\cap R$$ Denote by $T$ the group on the right. If $(R\cap M_k)N\supseteq T$ then $(R\cap M_k)(R\cap N)\supseteq T$ and by assumption $R\cap N\subseteq  M_k\cap N$ which implies that $M_k\supseteq T$, violating the fact that $M_1,\ldots, M_k$ are in general position. Hence $\pi(R\cap M_j)$ for $j>l$ are in general position.

\end{proof}

This lemma is very similar to Whiston's argument\footnote{Interested readers may refer to \cite{whiston2000maximal} or \cite{cameron2002independent} for more information}, but here we are working on a set of subgroups, hence readers may notice that the direction changes. Now we can prove Theorem \ref{THM1}:

\begin{proof}[Proof of Theorem \ref{THM1}] The proof is by the following strategy: given any set of maximal subgroups in general position, we replace pullback subgroups by standard subgroups in a way that does not change the total number of subgroups in question , so that the resulting set of maximal subgroups is still in general position. Let $\mathcal{S}=\{M_1,\ldots,M_k\}$ be a collection of maximal subgroups of $H\times K$ in general position, if no maximal subgroups are pullbacks then $k\leq \mathrm{MaxDim}(H)+\mathrm{MaxDim}(K)$ and we are done.  So suppose among $\mathcal{S}$ there is a pullback subgroup, say $M_k$ of the form $\Delta_{(N,N',\alpha)}$ induced by $N\lhd H$ and $N'\lhd K$ and $\alpha:H/N\cong K/N'\cong \mathbb{Z}_p$(By our assumption, pullback maximal subgroups must be of this form). Abbreviate it by $\Delta$, we know that necessarily $\Delta\cap (N\times N')=N\times N'$ and must contain the intersection of $M_j\cap (N\times N')$ for $j<k$. This, together with Lemma \ref{CruLem}, and the fact that $\mathbb{Z}_p^2$ is flat with $i=2$, indicates that we may assume without loss of generality, either $M_j\cap (N\times N')$ for $j<k$ are in general position, or $M_j\cap (N\times N')$ for $j<k-1$ are in general position and $M_{k-1}$ contains the intersection of them. Let us now consider separately these two cases.

\textit{Case 1}. Suppose $M_j\cap (N\times N')$ for $j<k$ are in general position. Let $$R=\bigcap\limits_{1\leq j<k}M_j$$ and $\pi: H\times K\to H\times K/(N\times N')$ be the projection, then $\pi(R)$ must be a nontrivial subgroup of $H\times K/(N\times N')\cong \mathbb{Z}_p^2$, otherwise $R\subseteq N\times N'\subseteq \Delta$, violating the fact that $\mathcal{S}$ is in general position. Because $\pi(R)$ is nontrivial, but $H/N\cap K/N'$ is trivial in the quotient, $\pi(R)$ is not contained in at least one of them, say $H/N$. We claim in this case that $\mathcal{S'}=\{M_1,\ldots, M_{k-1}, H\times N'\}$ is in general position. Indeed, since $M_i\cap (H\times N')$ for $1\leq i<k$ are in general position, and $H\times N'$ does not contain $R$ which is the intersection of $M_i$ for $1\leq i<k$, we see immediately that $\mathcal{S'}$ is in general position. 

\textit{Case 2}. Suppose $M_j\cap (N\times N')$ for $j<k-1$ are in general position and $M_{k-1}$ contains the intersection of them. Let $$R=\bigcap\limits_{1\leq j<k-1}M_j$$ and $\pi$ as in the first case, we claim that $\pi(R)=H/N\times K/N'$. Indeed, by Lemma \ref{CruLem}, $\pi(R\cap M_{k-1})$ and $\pi(R\cap\Delta)$ are in general position, which implies that $i(\pi(R))=2$. But the only subgroup of $\mathbb{Z}_p^2$ having $i=2$ is the whole group, hence the claim holds. Now it is not hard to show, by the same argument as in the first case, that the set $\mathcal{S'}=\{M_1,\ldots, M_{k-2}, H\times N', N\times K\}$ is in general position.

This shows that in any case we can replace pullback subgroups by standard subgroups, and a repeated application of this gives $k\leq \mathrm{MaxDim}(H)+\mathrm{MaxDim}(K)$. This is true for any collection of maximal subgroups of the direct product, and we conclude that $\mathrm{MaxDim}(H\times K)\leq \mathrm{MaxDim}(H)+\mathrm{MaxDim}(K)$. On the other hand, it is trivial that $\mathrm{MaxDim}(H\times K)\geq \mathrm{MaxDim}(H)+\mathrm{MaxDim}(K)$, hence $\mathrm{MaxDim}(H\times K)= \mathrm{MaxDim}(H)+\mathrm{MaxDim}(K)$. This completes the proof.

\end{proof}

Lemma \ref{CruLem} also yield the following:

\newtheorem{Flat}[LemofGou]{Proposition}

\begin{Flat}\label{Flat}
Suppose $H$ and $K$ are groups and $\mathrm{MaxDim}(K)=i(K)$, then:
\begin{center}
{ $\displaystyle \mathrm{MaxDim}(H\times K)=\mathrm{MaxDim}(H)+\mathrm{MaxDim}(K)$
}
\end{center}

\end{Flat}

\begin{proof}
Let $G=H\times K$ and let $M_1,\ldots, M_k$ be maximal subgroups that are in general position. Now after rearranging we may assume that $M_1\cap H,\ldots, M_l\cap H$ are in general position and $M_j$ contains their intersection for $j>l$. This together with Lemma 1.3 gives the bound $k-l\leq i(G/H)=i(K)=\mathrm{MaxDim}(K)$. Now for any maximal subgroup $M_i,~i\leq l$, since $\{M_i\cap H\vert1\leq  i\leq l\}$ is in general position, $M_i$ is either a standard subgroup induced by some maximal $M<H$, or a pullback. In both cases, $H\cap M_i$ will be an intersection of maximal subgroups. Indeed, if $M_i=\Delta_{(N,N',\alpha)}$, then $M_i\cap H=N$, which is the intersection of all maximal subgroups containing $N$ because $H/N$ is simple therefore Frattini Free; but if $M_i=M\times K$ then $M_i\cap H=M$. We show that $l\leq \mathrm{MaxDim}(H)$. Let $H_i=M_i\cap H$ for $1\leq i\leq l$. Now we claim that we can choose a maximal subgroup $M\supseteq H_1$ such that $M,H_2,\ldots,H_l$ are in general position. This is enough to show that $l\leq \mathrm{MaxDim}(H)$. Suppose for any such $M$, the sequence fail to be in general position, then the only possibility is that $M\supseteq \cap_{1<i\leq l}H_i$, since if $H_1\cap H_i$ are in general position for $1<i\leq l$, so are $M\cap H_i$ for any $M\supseteq H_1$. But this implies that $H_1\supseteq \cap_{1<i\leq l}H_i$ since $H_1$ equals the intersection of maximal subgroups containing it. This contradicts the fact that $H_i$ are in general position, so the claim holds. These arguments together give $k\leq \mathrm{MaxDim}(H)+\mathrm{MaxDim}(K)$. This is true for any set of maximal subgroups in general position, hence $ \mathrm{MaxDim}(G)\leq \mathrm{MaxDim}(H)+\mathrm{MaxDim}(K)$. Since the other inequality is trivial, the proof is complete.

\end{proof}

Theorem \ref{THM1} and Proposition \ref{Flat} cover a wide range of groups, e.g., all solvable groups, Mathieu groups, Janko groups $J_1$, $J_2$ and all symmetric and alternating groups satisfy assumptions of Theorem \ref{THM1} or Proposition \ref{Flat}\footnote{The results for Mathieu groups $M_{11}$ and $M_{12}$ is proved by T. Brooks(\cite{brooks2013generating}). The results for Mathieu group $M_{22}$ and Janko groups are joint work of Sophie Le, Tianyue Liu and me(to appear). The result for $M_{23}$ is from the work of Sophie Le(to appear), and the result for $M_{24}$ is from a remarkable theoretical proof of Tianyue Liu and R. Keith Dennis(\cite{MaxDimTianyue}). For alternating and symmetric groups, see \cite{whiston2000maximal} by J. Whiston.}, hence $\mathrm{MaxDim}(H\times K)=\mathrm{MaxDim}(H)+\mathrm{MaxDim}(K)$ if $K$ is one of the above groups. These all provide strong evidence that in general, $\mathrm{MaxDim}$ is additive on product of groups. 

It is not known yet in general, but in proving it or disproving it, one might need to consider $S\times S$ for $S$ a simple group. There are some reasons for this: first of all, $\mathrm{MaxDim}$ is additive automatically on $S\times S'$ if $S$ and $S'$ are simple but not isomorphic; second, for any group $G$, one has the conception of rumpf, denoted by $\mathfrak{R}(G)$, which is defined as the intersection of maximal normal subgroups, and it is easily seen that pullback subgroups are determined by $G/\mathfrak{R}(G)$, which is a product of simple groups. To this direction, the biggest family of non-flat simple groups is that of linear groups, and we deal with a part of this family in the following:

\newtheorem{SL2P}[LemofGou]{Proposition}

\begin{SL2P}\label{SL2P}
Let $\mathrm{PSL}(2,p)$ be the projective special linear group over the field with $p\geq 5$ elements, where $p$ is a prime. Then the following identity holds: $$\mathrm{MaxDim}(\mathrm{PSL}(2,p)\times \mathrm{PSL}(2,p))=2\mathrm{MaxDim}(\mathrm{PSL}(2,p))$$
\end{SL2P}

Before coming to the proof, it is worthwhile to introduce some important results about $\mathrm{PSL}(2,p)$ that will be used in the proof. The first one is Dickson's complete classification of maximal subgroups of $\mathrm{PSL}(2,p)$.

\newtheorem{Dickson}[LemofGou]{Theorem}

\begin{Dickson}[Dickson]\label{Dickson}
Maximal subgroups of $\mathrm{PSL}(2,p)$ comes from one of the following types:

(1). Point stabilizers $G_a$ for $a\in \mathbb{P}^1(\mathbb{F}_p)$, isomorphic to $\mathbb{Z}_p\rtimes \mathbb{Z}_{\frac{p-1}{2}}$.

(2). Stabilizers of a pair $G_{\{a,b\}}$ for $a,b\in \mathbb{P}^1(\mathbb{F}_p)$, isomorphic to $D_{p-1}$, dihedral group of order $p-1$.

(3). Stabilizers of a pair $G_{\{a,b\}}$ for $a,b\in \mathbb{P}^1(\mathbb{F}_{p^2})$, isomorphic to $D_{p+1}$, dihedral group of order $p+1$.

(4). Subgroups isomorphic to $A_4, S_4$ or $A_5$.
\end{Dickson}

Readers may refer to \cite{dickson2003linear} for more details. In \cite{whiston2002maximal}, J. Whiston and J. Saxl proved that $3\leq m(\mathrm{PSL}(2,p))\leq 4$. This is in fact also true for $\mathrm{MaxDim}$. In fact, in \cite{FinGrp}, D. Collins computed explicitly the intersection of maximal subgroups of the first three types. The result can be summarized as the following:

\newtheorem{Collins}[LemofGou]{Proposition}

\begin{Collins}[D. Collins]\label{Collins}
Intersections of maximal subgroups of $\mathrm{PSL}(2,p)$ of first three types are:

(1). $G_a\cap G_b\cong \mathbb{Z}_{(p-1)/2}$.

(2). $G_a\cap G_b\cap G_c$ is trivial.

(3). $G_a\cap G_{\{b,c\}}$ is either trivial or isomorphic to $\mathbb{Z}_2$.

(4). $G_{\{a,b\}}\cap G_{\{c,d\}}$ is isomorphic to $\mathbb{Z}_2\times \mathbb{Z}_2$, the Klein $4$-group.

\end{Collins}

In the following proof, we view $\mathrm{PSL}(2,p)$ both as matrices and as mobius functions on the space $\mathbb{F}_p\cup \{\infty\}$ interchangebly. In any case, $G_a$ will be point stabiliser, and $G_{\{a,b\}}$ will be the elements that stabilize the set $\{a,b\}$. 

\begin{proof}[Proof of Proposition \ref{SL2P}]
Denote by $G$ the group $\mathrm{PSL}(2,p)$. Let us first note that $\mathrm{MaxDim}(G)\leq 4$. Indeed, if $\mathrm{MaxDim}(G)>4$, then because $i(A_5)=i(S_4)=3$, the set of subgroups attaining maximal length would all come from first three types in Theorem \ref{Dickson}. But by Proposition \ref{Collins}, the intersection of any three maximal subgroups from first three types is a subgroup of $\mathbb{Z}_2$ if they are in general position. This is a contradiction. Hence $\mathrm{MaxDim}(G)\leq 4$.

We fix some notations $L_M=M\times G$ and $R_{M'}=G\times M'$, and $\Delta_\sigma=\{(g,\sigma(g)\mid g\in G)\}$. These are the only maximal subgroups of $G\times G$. Since the conclusion is trivial if no pullback subgroup is included, in what follows we assume that there is a pullback subgroup. Let $L_{M_1},\cdots,L_{M_k}$ and $R_{M_1'},\cdots,R_{M_s'}$ be standard maximal subgroups induced by $M_i,~M_j'<G$, and $\Delta_{id},\Delta_{\sigma_1},\cdots,\Delta_{\sigma_n}$ be pullback subgroups induced by $\mathrm{id},\sigma_i$, such that these subgroups are in general position(note that we can always assume that $\Delta_{\mathrm{id}}$ is included by applying appropriate automorphism), we wish to show that $k+s+n+1\leq 2\mathrm{MaxDim}(G)$. First note that $\{\Delta_{\mathrm{id}}\cap L_{M_i}, \Delta_{\mathrm{id}}\cap R_{M_j'}\}$ is in general position, then $\{M_i, M_j'\}$ is in general position in $\mathrm{PSL}(2,p)$. On the other hand, by Proposition \ref{Collins}, it is easily seen that if four maximal subgroups are in general position, then their intersection must be trivial. This indicates that if $k+s=4$, then $n=0$ and the cardinality of the whole set is less than $2\mathrm{MaxDim}(G)$.

So we may assume that $k+s$ does not attain $4$. But before working directly with all maximal subgroups, it seems convenient to first deal with $\Delta_{\sigma_i}$'s. In the following, we prove that in order that $\cap_{i}\Delta_{\sigma_i}\cap \Delta_{\mathrm{id}}$ be nontrivial, $n$ must be less than or equal to $4$.

We begin with pullback subgroups that are induced by inner homomorphism. Now let $\Delta_{id},\Delta_{h_1},\Delta_{h_2},\cdots,\Delta_{h_n}$ be maximal subgroups that are in general position, where $\Delta_{h_i}$ is induced by the conjugation of $h_i$, such that their intersection is not trivial. We show that $n\leq 3$. It is easily seen that their intersection is $\{(g,g)\mid g~\text{commutes with}~h_1,\cdots,h_k\}$. Also note that $\{h_1,\cdots,h_k\}$ is irredundant. Let $C(S)$ be the subgroup consists of elements commuting with $S$. Let $H$ be the subgroup generated by all $h_i$. If $H=G$ then since $G$ is simple, $C(G)$ is trivial and hence the intersection is trivial. This shows that in order that the subgroups above have nontrivial intersection, $H$ must be proper. We assume first that $H$ is contained in some maximal subgroup of the first three kind. But before going on, we will list here some observations that will be used throughout the proof:

(a). Suppose $h=\mathrm{diag}\{x,x^{-1}\}$, where $x^2\ne 1$, then $C(h)$ consists of only diagonal matrices iff $x^4\ne 1$. If $x^4=1$, then $C(h)$ consists of all diagonal matrices as well as matrices of the form $$\left(\begin{array}{cc} 0 & y\\ -y^{-1} & 0  \end{array}\right)$$ which we will denote by $\mathrm{Ad}\{y,-y^{-1}\}$.

(b). Suppose $h=\mathrm{Ad}\{y,-y^{-1}\}$ then $C(h)$ consists of matrices of the form $$\left(\begin{array}{cc} x & -zy^2\\ z & x  \end{array}\right)$$ as well as $$ \left(\begin{array}{cc} x & zy^2\\ z & -x  \end{array}\right)$$ where each matrix must have determinant one.

If $H<G_a$ for some $a$, then consider $C(h_1)$, it is a subgroup commuting with $h_1$. Since $h_1$ fixes $a$, it also fixes $ga$ for $g\in C(h_1)$, hence in order that $h_1$ be nontrivial, all elements in $C(h_1)$ can only possibly fix $a$ or move $a$ to another $b$(This is true by (2) of Proposition \ref{Collins}). If for all $h_i$, all elements in $C(h_i)$ fixes $a$, we show that $n\leq 2$. Under appropriate conjugation, we may assume that $a=\infty$ and $G_a \cong \mathbb{Z}_p\rtimes (\mathbb{Z}_p^\times)^2$ where the second factor acts by multiplying on the left(viewed as subgroup of the unit elements in $\mathbb{Z}_p$). Let $h_1=(x_1,y_1^2)$ and $h_2=(x_2,y_2^2)$. Let $(x,y^2)$ be an element that commute with both $h_1$ and $h_2$, then we have the formula:
$$\begin{cases} x_1+y_1^2x=x+y^2x_1\\
x_2+y_2^2x=x+y^2x_2 \end{cases}$$
and more precisely the following matrix equation:
$$\left( \begin{array}{cc} 1-y_1^2 & x_1\\ 1-y_2^2 & x_2 \end{array} \right)
\left( \begin{array}{c} x \\ y^2 \end{array} \right)=\left( \begin{array}{c} x_1 \\ x_2 \end{array} \right)
$$

If the matrix on the left has determinant zero, then the two rows are linearly dependent, in which case $C(h_1)= C(h_2)$, which cannot happen since the subgroups are in general position. Hence the only possibility is that the matrix is invertible, which means that there is only one possible solution for the above equation, i.e., $(0,1)$. Hence $C(h_1)\cap C(h_2)=1$ and $n\leq 2$.

If $C(h_1)$ consists of elements fixing $a$ and switching $\{a,b\}$, then $h_1\in G_a\cap G_b$, and under appropriate conjugation, we may assume that $h_1$ is diagonal. If $h_2$ is not in this intersection, then $h_2$ maps $b$ to some $c$, in which case any element $g$ in $C(h_1)\cap C(h_2)$ will need to fix $\{a,b\}$ and satisfy $gc=gh_2 b=h_2 gb $ which equals either $c$ or $a$. But since $g$ fixes $\{a,b\}$, $g$ must also fix $c$, in which case $C(h_1)\cap C(h_2)\subseteq \mathbb{Z}_2$ and nothing else can be added otherwise the intersection will be trivial. Now consider the other case, where all $h_i$ lies in $G_a\cap G_b$, and hence all $h_i$ are diagonalisable. In this case, by the observation (1), we conclude that $C(h_1)=C(h_2)$ or $C(h_1)\supseteq C(h_2)$ or $C(h_1)\subseteq C(h_2)$. This contradicts the fact that the subgroups are in general position.

Up to now the first case is all done and the conclusion is that whenever $H<G_a$, $n\leq 2$, or to say, for any three of the $C(h_i)$, if they are in general position, their intersection must be trivial.

Now consider the case where $H<G_{\{a,b\}}$. Here we may well work at $\mathrm{PSL}(2,\mathbb{F}_{p^2}\!)$ which unifies type (2) and (3) in Theorem \ref{Dickson}, because the above observations are still true here. Without loss of generality, we may assume that $H$ consists of matrix of the form $\mathrm{diag}\{x,x^{-1}\}$ and $\mathrm{Ad}\{y,-y^{-1}\}$. By same argument as in the above case, we know that there will at most be one $h_i$ of diagonal form. Suppose $h_1$ is of diagonal form $\mathrm{diag}\{x_1,x_1^{-1}\}$ and $h_2$ is of the form $\mathrm{Ad}\{y_2,-y_2^{-1}\}$. Direct computation by observation (2) shows that $C(h_1)\cap C(h_2)$ either contains only $\mathrm{diag}\{x,x^{-1}\}$ with $x^4=1$(if $C(h_1)$ contains only diagonal matrices) or contains these together with those $\mathrm{Ad}\{y,-y^{-1}\}$ such that $y^4=y_2^4$. In the first case, $C(h_1)\cap C(h_2)=\mathbb{Z}_2$ and nothing more can be added to the sequence. In the second case, if we choose $ h_3=\mathrm{Ad}\{y_3,-y_3^{-1}\}$ with $y_3^4=y_2^4$, then $C(h_2)\cap C(h_1)=C(h_3)\cap C(h_1)$ which contradicts the fact that they are in general position. Hence we must choose $h_3$ such that $y_3^4\ne y_2^4$, which means that $C(h_1)\cap C(h_2)\cap C(h_3)$ contains only elements of the form $\mathrm{diag}\{x,x^{-1}\}$ with $x^4=1$, and so $C(h_1)\cap C(h_2)\cap C(h_3)=\mathbb{Z}_2$. And we cannot add any more $h_i$ in otherwise the intersection will be trivial. Finally, if we only choose element of anti-diagonal form, and assume that we have chosen $h_i=\mathrm{Ad}\{y_i,-y_i^{-1}\}$ for $1\leq i\leq 3$, and $y_i^4$ are not all equal, then it is possible using same argument as above to show that $C(h_1)\cap C(h_2)\cap C(h_3)\subseteq \mathbb{Z}_2$. This finishes the case where $H<G_{\{a,b\}}$.

If $H<S_4$ or $A_5$, and if $n=3$, then the $h_i$'s form an irredundant generating sequence for respective groups. Because in either case the group is strongly flat and has trivial center, the intersection will have to be trivial. Hence in order that the intersection being nontrivial, $n\leq 2$. The same is true if $H<A_4$.

The conclusion here is that if $\Delta_{\mathrm{id}}, \Delta_{h_1}, \Delta_{h_2},\Delta_{h_3}$ are in general position, there intersection is a subset of $\mathbb{Z}_2$. We now apply an argument similar to Whiston's to pass the information to $\mathrm{Aut}(G)$\footnote{Here we use the fact that $\mathrm{Out}(G)=\mathbb{Z}_2$.}. Suppose we have $\Delta_{\sigma_1},\cdots,\Delta_{\sigma_4}$ together with $\Delta_{\mathrm{id}}$ are in general position, and one of the $\sigma_i$, say $\sigma_1$ is not an inner-morphism. Then we can multiply $\sigma_1$ appropriately to $\sigma_j$'s to make them inner, say $\sigma_1^{s_i}\sigma_j$ is inner, and we have $\mathrm{Rad}(\Delta_{\sigma_1},\ldots,\Delta_{\sigma_4},\Delta_{\mathrm{id}})\footnote{$\mathrm{Rad}$ means the radical, i.e., taking intersection.}=\mathrm{Rad}(\Delta_{\sigma_1},\Delta_{\sigma_1^{s_2}\sigma_2},\cdots,\Delta_{\sigma_1^{s_4}\sigma_4},\Delta_{\mathrm{id}})$. Consider $\Delta_{\mathrm{id}}\cap \Delta_{\sigma_1^{s_j}\sigma_j}~\text{for}~j>1$. We claim that they are in general position. If not, say $\Delta_{\mathrm{id}}\cap \Delta_{\sigma_1^{s_1}\sigma_2}\cap \Delta_{\sigma_1^{s_3}\sigma_3}<\Delta_{\mathrm{id}}\cap \Delta_{\sigma_1^{s_4}\sigma_4}$, then $\mathrm{Rad}(\Delta_{\sigma_1},\Delta_{\sigma_1^{s_2}\sigma_2},\cdots,\Delta_{\sigma_1^{s_4}\sigma_4},\Delta_{\mathrm{id}})=\mathrm{Rad}(\Delta_{\sigma_1},\Delta_{\sigma_1^{s_2}\sigma_2},\Delta_{\sigma_1^{s_3}\sigma_3},\Delta_{\mathrm{id}})$, which implies $\cap_i \Delta_{\sigma_i}\cap \Delta_{\mathrm{id}}=\cap_{i<4} \Delta_{\sigma_i}\cap \Delta_{\mathrm{id}}$. This contradicts the fact that the five subgroups are in general position. But if $\Delta_{\mathrm{id}}\cap \Delta_{\sigma_1^{s_j}\sigma_j}~\text{for}~j>1$ are in general position, then their intersection lies inside $\mathbb{Z}_2$ as has been shown above. Hence if $\Delta_{\sigma_i}$ for $1\leq i\leq 4$ and $\Delta_{\mathrm{id}}$ are in general position, their intersection is a subset of $\mathbb{Z}_2$.

Now we can finish the proof of the proposition. Suppose we are given some collections of maximal subgroups $L_{M_1},\cdots,L_{M_k}$ and $R_{M_1'},\cdots,R_{M_s}'$ with $k+s\leq 3$, $\Delta_{\mathrm{id}}$, $\Delta_{\sigma_1},~\Delta_{\sigma_2},~\Delta_{\sigma_3},~\Delta_{\sigma_4}$ where $\sigma_i$'s are nontrivial automorphisms. We are done if $\mathrm{MaxDim}(G)=4$ since the argument above indicates that any time if pullback subgroups are included, the total number of maximal subgroups has to be no more than $8$. If $\mathrm{MaxDim}(G)=3$, then in order that the equality fails, all but one of the maximal subgroups of the above need to be included. We show that this is impossible. Suppose we use all the $\Delta_{\sigma_i}$, then the above argument shows that the intersection of these subgroups with $\Delta_{\mathrm{id}}$ is either $\mathbb{Z}_2$ or trivial, indicating that we can only possibly add one more in the set. So no more than three $\Delta_{\sigma_i}$'s other than $\Delta_{\mathrm{id}}$ can be included, and in order to attain $7$, we need to use all  standard maximal subgroups above, that is $L_{M_i}$ and $R_{M_j'}$. If the $M_i, M_j'$ are all from first three types, then by easy computation, their intersection will be either trivial or isomorphic to $\mathbb{Z}_2$, in which case only $5$ maximal subgroups can be included. Hence we must have some $A_5$ or $S_4$ in $M_i, M_j'$. But in this case, since $i(A_5)=i(S_4)=3$ we conclude that we can possibly have $3$ more choices of other types other than $\Delta_{\mathrm{id}}$. In either case, $k+s+n+1\leq 6\leq 2\mathrm{MaxDim}(G)$. This is true for any set of maximal subgroups in which there is one pullback subgroup, and because it is trivially true if no pullback subgroup is used, we conclude that $\mathrm{MaxDim}(G\times G)\leq 2\mathrm{MaxDim}(G)$. Since the other inequality is trivial, the proof is complete.

\end{proof}

We conclude this section by some comments. First of all, it seems to the author that in proving or disproving $\mathrm{MaxDim}(S^2)=2\mathrm{MaxDim}(S)$ for simple $S$ one will inevitably use classification theorem. Second, it is not hard to show that for a simple $S$, if $\mathrm{MaxDim}(S^2)=2\mathrm{MaxDim}(S)$, then $\mathrm{MaxDim}(S^n)=n\mathrm{MaxDim}(S)$. Third, even if the identity holds for simple groups, there is no obvious evidence that it should be true in general. So before going into the proof for all simple groups, I think it is worthwhile to build the bridge from simple groups to general groups. R. Keith Dennis, on the other hand, suggested that one should first focus on simple groups and look for stronger properties than that in Proposition \ref{SL2P}, among which there might be illustrations toward general situations. This is possibly doable because in Proposition \ref{SL2P}, it seems very likely that once a pullback subgroup is included, the number of subgroups in general position is no greater than $\mathrm{MaxDim}(S)+1$. It will be great if such conclusion does hold, but by now we only know this holds for a small number of groups, i.e., flat simple groups. This is the following proposition proved by R. Fernando in Cornell 2017 summer REU program:

\newtheorem{RavProp}[LemofGou]{Proposition}

\begin{RavProp}[R. Fernando]
Let $S$ be a simple group. If $M_1,\ldots, M_k$ are maximal subgroups of $S\times S$ that are in general position, then we have:
$$k\leq \begin{cases}
~~2\mathrm{MaxDim}(S)\hspace{2cm}\text{if $M_i$ are all standard}\\
\\
~~i(S)+1\hspace{2.75cm}\text{if some $M_i$ is a pullback}
\end{cases}
$$
\end{RavProp}
\begin{proof}
If all $M_i$ are of standard form, it is trivial that one has $k\leq 2\mathrm{MaxDim}(S)$. If some $M_i$, say $M_1$ is a pullback, then $M_1\cong S$ and because $M_i\cap M_1$ are in general position for $i>1$, we have $k-1\leq i(S)$, hence $k\leq i(S)+1$. This completes the proof.

\end{proof}

\section{Relative versions of $\mathrm{MaxDim}$}\label{3}

This section serves as a transition between Section \ref{2} and Section \ref{4}, some definitions and results will be introduced and will be used in next section. In order to further understand Maximal Dimension, inspired by \cite{genseqfingrp}, we come to the definition of relative versions of this function. The first one is the following:

\newtheorem{DefRel1}{Definition}[section]
\begin{DefRel1}\label{DefRel1}
Let $G$ be a group and $H<G$, then the \textbf{maximal dimension of $G$ relative to $H$}, which we denote by $\mathrm{MaxDim}(G,H)$, is defined by:
$$\displaystyle \mathrm{MaxDim}(G,H)=\max_\mathcal{S} |\mathcal{S}|$$
where the maximum is taken over all $\mathcal{S}=\{M_1,\ldots,M_k\}$, collections of maximal subgroups of $G$ such that $\{M_1\cap H,\ldots,M_k\cap H\}$ is in general position.
\end{DefRel1}

The naturality of this definition can be seen in Lemma \ref{CruLem}. There is another version, which comes naturally when studying group extensions:
\newtheorem{DefRel2}[DefRel1]{Definition}
\begin{DefRel2}
Let $G$ be a group and $H$ another group acting on $G$. Define the \textbf{maximal dimension of $G$ under $H$}, which is denoted by $\mathrm{MaxDim}_H(G)$, by:
$$\displaystyle \mathrm{MaxDim}_H(G)=\max_\mathcal{S} |\mathcal{S}|$$
where the maximum is taken over all $\mathcal{S}=\{M_1,\ldots,M_k\}$, collections of maximal $H$-invariant subgroups of $G$ that are in general position.
\end{DefRel2}

Our goal in this and next section is to use these two conceptions to study behaviors of maximal dimension under group extensions, especially extensions by abelian groups. To this vein, we prove the following:

\newtheorem{PropRel}[DefRel1]{Proposition}
\begin{PropRel}\label{PropRel}
Let $N\lhd H$ and $\pi:H\to H/N$. Let $M_1,\ldots,M_k$ be maximal subgroups of $H$ that are in general position, such that $M_1\cap N,\ldots,M_l\cap N$ are in general position and $M_j$ contains $\cap_{1\leq i\leq l}M_i\cap N$ for any $j>l$. Let $R=\cap_{1\leq i \leq l}M_i$ and $\pi:H\to H/N$ be the projection. Then $l\leq \mathrm{MaxDim}(H)-\mathrm{MaxDim}(H/N,\pi(R))$. Moreover, if $N/\Phi(N)$ is abelian, $l\leq \mathrm{MaxDim}(H)-\mathrm{MaxDim}(H/N)$, and $l\leq \mathrm{MaxDim}_H(N)$.
\end{PropRel}

\begin{proof}
For the first statement, suppose $M_1',\ldots,M_s'$ are maximal subgroups of $H$ containing $N$ such that $\pi(M_1')\cap \pi(R),\ldots,\pi(M_s')\cap \pi(R)$ are in general position, it is not hard to show that $\{M_1,\ldots,M_l,M_1',\ldots,M_s'\}$ are in general position, hence $l\leq \mathrm{MaxDim}(H)-s$. Taking $s$ to be maximal, we get $$l\leq \mathrm{MaxDim}(H)-\mathrm{MaxDim}(H/N,\pi(R))$$ 

For the second statement, since $\Phi(N)<\Phi(H)$, we assume that $H$ is Frattini free, which implies by our assumption that $N$ is abelian. Note that if $M$ is a maximal subgroup of $H$ not containing $N$, we claim that $M\cap N$ is a maximal $H$-invariant subgroup of $N$. It is indeed $H$-invariant because $MN=H$, so we only need to show that $M\cap N$ is maximal. If not, say $ M\cap N\subsetneq N'\subsetneq N$ for some $N'\lhd H$. Since $M$ is maximal, $MN'=H$ and $H/N=M/(M\cap N)=MN'/N'=H/N'$, a contradiction. Hence $M\cap N$ is maximal $H$-invariant. This gives the bound $l\leq \mathrm{MaxDim}_H(N)$. Now we proceed to prove that $\pi(R)=H/N$, which will complete the proof of the Proposition. Suppose $R'$ is a subgroup of $H$ such that $\pi(R')=H/N$, and $M$ a maximal subgroup of $H$ not containing $N$ such that $M\cap N\nsupseteq R'\cap N$. We claim that $\pi(M\cap R')=H/N$. First of all, since $M\cap N$ is maximal invariant, $(M\cap N)(R'\cap N)=N$. Now for any $h\in H/N$, there exists $m\in M$ and $r\in R'$ such that $\pi(m)=\pi(r)=h$. This indicates that $mr^{-1}\in N=(M\cap N)(R'\cap N)$, hence there exists $n_1\in M\cap N$ and $n_2\in R'\cap N$ such that $mr^{-1}=n_1^{-1}n_2$, which indicates that $n_1m=n_2r\in R'\cap M$. This element has image $h$. This shows that every element in $H/N$ has a pre-image in $R'\cap M$, hence $\pi(R'\cap M)=H/N$, and the proof is complete.  

\end{proof}

Proposition \ref{PropRel} gives several bounds on $l$, a value coming from Lemma \ref{CruLem} depending on the collection of maximal subgroups. There is a corollary to this which is also related to the previous section:
\newtheorem{CorAbRf}[DefRel1]{Corollary}
\begin{CorAbRf}
Suppose $H$ and $K$ are groups with $\mathfrak{R}(H)$ and $\mathfrak{R}(K)$ abelian, and both $H/\mathfrak{R}(H)$ and $K/\mathfrak{R}(K)$ are product of flat simple groups. Then $\mathrm{MaxDim}(H\times K)=\mathrm{MaxDim}(H)+\mathrm{MaxDim}(K)$.

\end{CorAbRf}

\begin{proof}
Let $M_1,\ldots,M_n$ be maximal subgroups of $H\times K$ that are in general position, and consider $N=\mathfrak{R}(H)\times \mathfrak{R}(K)$. Suppose $M_1\cap N,\ldots,M_k\cap N$ are in general position, and $M_j$ contains the intersection of them. These $M_i$ must be type (1) subgroups either not containing $\mathfrak{R}(H)\times K$(if they are induced by $M<H$), or not containing $H\times\mathfrak{R}(K)$(if they are induced by $M'<K$). Hence by Proposition 3.5, $k\leq \mathrm{MaxDim}(H)+\mathrm{MaxDim}(K)-\mathrm{MaxDim}(H/\mathfrak{R}(H))-\mathrm{MaxDim}(K/\mathfrak{R}(K))$. On the other hand, by Lemma 1.3, $n-k\leq i(H/\mathfrak{R}(H)\times K/\mathfrak{R}(K))=i(H/\mathfrak{R}(H))+i(K/\mathfrak{R}(K))$, which is equal to $\mathrm{MaxDim}(H/\mathfrak{R}(H))+\mathrm{MaxDim}(K/\mathfrak{R}(K))$ by flatness. Combining these inequalities we get the desired result.

\end{proof}

Proposition \ref{PropRel} also yield the following corollary, which will be used in next section.

\newtheorem{Corsplit}[DefRel1]{Corollary}

\begin{Corsplit}\label{Corsplit}
Suppose $H=N\rtimes T$ where $N$ is abelian and $T$ is flat. Then $\mathrm{MaxDim}(H)=\mathrm{MaxDim}_T(N)+\mathrm{MaxDim}(T)$

\end{Corsplit} 

\begin{proof}
Note that Proposition 3.5 gives the bound $\mathrm{MaxDim}(H)\leq \mathrm{MaxDim}_T(N)+\mathrm{MaxDim}(T)$. But since the sequence splits, it is very easy to demonstrate a sequence of that length, which consist of $N_1\rtimes T,\ldots,N_k\rtimes T$ and $N\rtimes M_1,\ldots,N\rtimes M_t$ where $k=\mathrm{MaxDim}_T(N)$ and $t=\mathrm{MaxDim}(T)$.

\end{proof}

\section{Groups with $\mathrm{MaxDim}=m$ and groups with \\$\mathrm{MaxDim}>m$}\label{4}

In this section, we will give another proof, using the results and techniques in Section \ref{3}, the theorem that has been proved by E. Detomi and A. Lucchini in \cite{detomi2016maximal}, stating that if $G'$ is nilpotent, then $m(G)=\mathrm{MaxDim}(G)$.  Also, we will prove a proposition conjectured by Ellie Thieu in 2017 Cornell summer SPUR forum, which not only exhibits a connection between groups with $m<\mathrm{MaxDim}$ and groups with $m<i$, but also explains the reason that the group $R=\mathrm{smallgroup}(720,774)$ found by R. Fernando in \cite{fernando2015inequality} is the smallest solvable group with $m<\mathrm{MaxDim}$. In fact, the proof of this proposition is based on computations on GAP for $R$.

To begin with, we introduce the following version of Goursat lemma:

\newtheorem{LemGst}{Lemma}[section]
\begin{LemGst}[Goursat Lemma for $G$-groups]\label{LemGst}
Let $G$ be a group acting on $H$ and $K$. Let $\pi_1:H\times K\to H$ and $\pi_2:H\times K\to K$. Suppose $\Delta$ is a proper $G$-invariant subgroup of $H\times K$ such that $\pi_1(\Delta)=H$ and $\pi_2(\Delta)=K$, then there exists invariant normal subgroups $N\lhd H$ and $N'\lhd K$, and a $G$-isomorphism $\alpha:H/N\to K/N'$ such that $\Delta=\{(h,k)\vert \alpha(\bar{h})=\bar{k}\}$.  
\end{LemGst}

\begin{proof}
Let $N=\Delta\cap H\times 1$. This is a $G$-invariant subgroup, and by definition a normal subgroup of $\Delta$, being the kernel of the projection $\pi_2$. Since $\Delta$ projects on the first factor, $N\lhd H$. The same reason indicates that $N'=\Delta\cap 1\times K$ is normal. The rest of the statements will be routine check and will be omitted.
\end{proof}

Another lemma important to us is a structure theorem for groups with nilpotent commutator subgroups:
\newtheorem{LemSup}[LemGst]{Lemma}
\begin{LemSup}\label{LemSup} 
Let $G$ be a finite group with $G'$ nilpotent, and suppose $\Phi(G)=1$. Then there are abelian subgroups $F,~H<G$ such that $F\lhd G$ and $G=F\rtimes H$. Moreover, $F$ can be decomposed into minimal $H$-groups. 
\end{LemSup}

\begin{proof}
See \cite{detomi2016maximal}.

\end{proof}

\newtheorem{ThmSup}[LemGst]{Theorem}
\begin{ThmSup}
Let $G$ be a group with $G'$ nilpotent, then $\mathrm{MaxDim}(G)=m(G)$.
\end{ThmSup}

\begin{proof}

Write $G=F\rtimes H$ as in Lemma \ref{LemSup} and assume that $G$ is Frattini Free, which implies that $F$ can be decomposed into minimal $H$ vector spaces. Since $H$ is abelian, Corollary \ref{Corsplit} is applied and we get that $\mathrm{MaxDim}(G)=\mathrm{MaxDim}_H(F)+m(H)$. Since $m(G)=m(H)+\text{number of irreducible $H$ components of $F$}$\footnote{See, for example, \cite{lucchini2013largest}, or \cite[p. 52]{genseqfingrp} Proposition 3.4.1},  we are only left to prove that $\mathrm{MaxDim}_H(F)$ is exactly the number of complemented chief factors contained in $F$(i.e., the number of irreducible components of $F$). For this purpose, we may assume without loss of generality, that $F=V^n$ for some minimal $H$ module $V$ and some $n>0$. 

Let $U$ be a maximal $H$-invariant subgroup of $F$, then by Lemma \ref{LemGst}, $U$ either discard one of the factors $V$ totally(i.e., $U=V^{n-1}$), or is determined by diagonally linking two copies of them(by some automorphism of $V$). We show by induction that the number of maximal invariant subgroups used should always be less than or equal to the number of irreducible component. Now take any maximal subgroup $U$ that is used, we deviate two cases. First, suppose $U=V^{n-1}$, i.e., $U$ is of normal type. For any other $W$, since $W\ne U$, $W$ either discard another copy of $V$, or is determined by diagonally linking two copies of $V$. In the first case, $W\cap U$ is clearly maximal in $U$. In the second case, if $W$ links two copies of $V$ that are both inside of $U$, then $W\cap U$ is a maximal invariant subspace of $U$ determined by diagonally linking two copies of $V$. Otherwise $W$ links two copies of $V$ one of which is the complement of $U$. In this case, $W\cap U$ is the direct sum of $V^{n-1}$ and the subspace of $V$ consisting of fixed points of a nontrivial $H$-automorphism. Since $V$ is minimal this latter space is trivial. Hence $W\cap U$ is maximal again. In either case $W\cap U$ is maximal and inductive hypothesis can be applied to give the bound. Second, suppose $U$ is a diagonal subspace. In this case we also have that $U\cong V^{n-1}$. For any other $W$, if $W$ is a normal form or is determined by diagonally linking a pair different from that of $U$, then $W\cap U$ is maximal. If $W$ is determined by diagonally linking the same pair as $U$, but in a different way, then $W\cap U$ is again the direct sum of $V^{n-2}$ with a diagonal subgroup determined by the fixed point space of a nontrivial $H$-automorphism of $V$, which is trivial by the same reason as above. Hence $W\cap U$ is maximal, and the inductive hypothesis can be applied. This proves that $\mathrm{MaxDim}_H(F)$ is indeed the number of irreducible components of $F$, and the proof is complete.
  
\end{proof}

In \cite{genseqfingrp}, D. Collins used relative versions of $m$ to deal with $m(G)-m(G/N)$, by which he proved that $m(G\times H)=m(G)+m(H)$. In \cite{lucchini2013largest}, A. Lucchini investigated more closely the quantity $m(G)-m(G/N)$ for minimal normal $N$, and showed that this quantity only depends on the action of $G$ on $N$. The above results come from the same idea, i.e., to look at $\mathrm{MaxDim}(G)-\mathrm{MaxDim}(G/N)$. But A. Lucchini proved(in private communication) that this value can be arbitrarily large even for minimal normal abelian $N$, so we are still not close to fully understanding this value. Next result is based on computations for $R$, which also gives a way to construct infinite class of groups with $m<\mathrm{MaxDim}$. 

\newtheorem{PropEl}[LemGst]{Proposition}
\begin{PropEl}\label{PropEl}
Let $S$ be a non flat group. There exists a semi-simple $S$-group $V$ such that $\mathrm{MaxDim}(V\rtimes S)=\mathrm{MaxDim}_S(V)+i(S)$.

\end{PropEl}

Clearly because $m(V\rtimes S)=\mathrm{MaxDim}_S(V)+m(S)$, if $S$ is flat, this indeed gives an example of groups with $m<\mathrm{MaxDim}$, and somehow they are as many as groups with $m<i$. We need a lemma to prove the result:

\newtheorem{LemEl}[LemGst]{Lemma}
\begin{LemEl}\label{LemEl}
Let $H<S$ be a subgroup, and $p$ a prime not dividing the cardinality of $S$. Then there exists a finite $S$-space $V$ over $\mathbb{F}_p$ and a vector $v\in V$ such that $\mathrm{Stab}(v)=H$.

\end{LemEl}

\begin{proof}
Let $V_0=\mathbb{F}_p$, viewed as the trivial $H$ module. Let us consider the induced module $V=\mathrm{ind}_H^S(V_0)$. An element $f\in V$ is defined as a map $f:S\to V_0$ such that $f(gh)=h^{-1}(f(g))=f(g)$ since $V_0$ is a trivial $H$ module, which is equivalent to that $f$ takes constant value along the left cosets of $H$. Now let us choose such a function $f\in V$, taking $0$ on $H$ and $1$ on all other left cosets. This is a nontrivial element in $V$. We claim that $\mathrm{Stab}(f)=H$. Pick $s\in \mathrm{Stab}(f)$, then by definition $sf(g)=f(g)$ for all $g\in S$. Taking $g=1$ we see that $sf(1)=f(s^{-1})=f(1)$, hence $f(s^{-1})=0$ and so $s\in H$. On the other hand, for any $h\in H$, $hf(g)=f(h^{-1}g)=f(h^{-1}gh)=f(g)$ as can be easily seen. Hence indeed $\mathrm{Stab}(f)=H$ and the proof is complete.

\end{proof}

\begin{proof}(of Proposition \ref{PropEl})
Fix a prime $p$ not dividing the cardinality of $S$, and let $H_1,\ldots,H_n$ be subgroups of $S$ that are in general position with $n=i(S)$. For each $i$, there exists $V_i$ over $\mathbb{F}_p$ and $v_i\in V_i$ such that $\mathrm{Stab}(v_i)=H_i$. For each $i$, we can decompose $V_i$ into irreducible subspaces and the elements that stabilize $v_i$ stabilize all factors in the decomposition. Meanwhile, if $\cap_{1\leq i\leq t} M_i, N_2,\ldots, N_k$ are in general position, there must be some $i$ such that $M_i, N_2,\ldots, N_k$ are also in general position. Hence we may assume that each $v_i$ spans an irreducible component of $V_i$ under the action of $S$. Let $V=V_1\oplus V_2\oplus\cdots\oplus V_n$, and consider $G=V\rtimes S$. First of all, it is easily seen using Lemma 1.3 that $\mathrm{MaxDim}(G)\leq \mathrm{MaxDim}_S(V)+i(S)$, so we must prove the other direction. We do this by exhibiting a sequence of this length. 

Fix a decomposition of $V$, say $W_1\oplus W_2\oplus\cdots \oplus W_m$ such that the first $n$ factors are just $W_i=\mathrm{span}_S(v_i)$, the irreducible subspaces generated by each $v_i$. Define maximal subgroups of $G$ by $$M_i=(\bigoplus\limits_{j\ne i}W_j)\rtimes S$$ for $1\leq i\leq m$ and $$T_k=(\bigoplus\limits_{j\ne k}W_j) v_kSv_k^{-1}$$ for $1\leq k\leq n$. We claim that $M_i, T_k$ are in general position. Let $\pi:G \to S$ be the projection. First, any $T_k$ does not contain the intersection of the other, since $\pi(\cap_{1\leq i\leq n}M_i\cap T_k)=\pi(\{v_ksv_k^{-1}\vert s\in H_k\})=H_k$ and we know that all $H_k$ are in general position. Next, $M_j$ for $j>n$ does not contain the intersection of the other, since the intersection of all except $M_j$ contains $W_j$ but $M_j$ doesn't. Finally, $M_i$ for each $1\leq i\leq n$ does not contain the intersection of the other. Without loss of generality, consider $M_1$, and let $R=\cap_{j>1}M_j\cap T_1$, by calculation this is just $v_1Sv_1^{-1}$, and we see that $\pi(R\cap M_1)=H_1$. What is $\pi(R\cap T_k)$ for $k>1$? For instance, for any $v_1sv_1^{-1}\in R$, if it is in $T_k$, then $(v_1sv_1^{-1}s^{-1})s=(wv_ksv_k^{-1}s^{-1})s$ for some $w\in \oplus_{j\ne k}W_j$, and so $v_1sv_1^{-1}s^{-1}-v_ksv_k^{-1}s^{-1}=w$ for some $w\in \oplus_{j\ne k}W_j$. But since the $k$-th component of $w$ is trivial by definition, $v_ksv_k^{-1}s^{-1}=0$ and so $s\in H_k$. Now it is easily seen that $\pi(R\cap T_k)=H_k$ for $k>1$, and because $H_k$ are in general position, $M_1\cap R$ does not contain the intersection of $R\cap T_k$ for $k>1$. This indicates that the subgroups are in general position. Counting the number of subgroups in question we see that $\mathrm{MaxDim}(V\rtimes S)\geq \mathrm{MaxDim}_S(V)+i(S)$, and the result follows from these two inequalities.

\end{proof}

In fact, the pattern of subgroups of $R$ is exactly of this form, and the reason $m(R)<\mathrm{MaxDim}(R)$ is because $R$ has a non flat top $S$(smallest non flat group, in fact, by computation in GAP) and a bottom $V$ which $S$ acts on, such that stabilizers of elements in $V$ play a role in making $\mathrm{MaxDim}$ one larger. This also shows that in Lemma \ref{CruLem}, the subgroups $\pi(R\cap M_j)$ are not predictable in the quotient, in particular not necessarily maximal.

\newpage
\section{Acknowledgement}

The author wants to take this chance to express his gratitude for many people, especially for Professor R. Keith Dennis and graduate mentor Ravi Fernando, for their encouragements and advise, as well as all Cornell 2017 REU participants, for many insightful comments and fruitful discussions. The author also wants to thank Professor Andrea Lucchini for his private communication with the author, giving the author many new perspectives on this topic. The author is also very thankful for Michael Vincent Crea for being a friend to the author and providing much support during this summer.

\newpage

\bibliographystyle{plain}

\bibliography{Final}

\end{document}